\documentclass[12pt]{article}
\usepackage[a4paper,margin=0.5in]{geometry}
\usepackage[titletoc,title]{appendix}
\usepackage{changepage}
\usepackage[utf8x]{inputenc}
\usepackage{textcomp,marvosym}
\usepackage{fixltx2e}
\usepackage{amsmath,amssymb}
\usepackage{cite}
\usepackage{nameref,hyperref}
\usepackage{authblk}
\usepackage{algorithm}
\usepackage{algpseudocode}
\usepackage{dsfont}
\usepackage{xcolor}
\usepackage{enumitem}
\usepackage{caption,setspace}
\usepackage{subcaption}
\usepackage{pdflscape}
\usepackage{graphicx}
\usepackage{placeins}
\usepackage{mathtools}

\usepackage{tikz}
\usepackage{pgfplots}

\renewcommand{\eqref}[1]{Equation~(\ref{#1})}

\hypersetup{pdfpagemode=UseNone}

\title{{\sc Travelling waves, blow-up and extinction \\ in the Fisher-Stefan model}}

\author[1]{Scott~W. McCue\footnote{To whom correspondence should be addressed. E-mail: scott.mccue@qut.edu.au}}
\author[1]{Maud~El-Hachem}
\author[1]{Matthew~J. Simpson}

\affil[1]{School of Mathematical Sciences, Queensland University of Technology \newline Brisbane QLD 4001, Australia}

\begin{document}

\maketitle
\begin{abstract}
While there is a long history of employing moving boundary problems in physics, in particular via Stefan problems for heat conduction accompanied by a change of phase, more recently such approaches have been adapted to study biological invasion.  For example, when a logistic growth term is added to the governing partial differential equation in a Stefan problem, one arrives at the \textit{Fisher-Stefan} model, a generalisation of the well-known Fisher-KPP model, characterised by a leakage coefficient $\kappa$ which relates the speed of the moving boundary to the flux of population there.  This Fisher-Stefan model overcomes one of the well-known limitations of the Fisher-KPP model, since time-dependent solutions of the Fisher-Stefan model involve a well-defined front which is more natural in terms of mathematical modelling.  Almost all of the existing analysis of the standard Fisher-Stefan model involves setting $\kappa > 0$, which can lead to either invading travelling wave solutions or complete extinction of the population.  Here, we demonstrate how setting $\kappa < 0$ leads to retreating travelling waves and an interesting transition to finite-time blow-up.  For certain initial conditions, population extinction is also observed.  Our approach involves studying time-dependent solutions of the governing equations, phase plane and asymptotic analysis, leading to new insight into the possibilities of travelling waves, blow-up and extinction for this moving boundary problem.  Matlab software used to generate the results in this work are available on  \href{https://github.com/ProfMJSimpson/Fisher_blowup}{Github}.
\end{abstract}

\noindent
Keywords: Fisher-KPP equation; Stefan problem; moving boundary problem; retreating fronts; sharp-fronted travelling waves; blow-up; extinction

\newpage
\section{Introduction} \label{intro}

Moving boundary problems are frequently used to model phenomena in physics, in particular melting/freezing processes where the moving boundary represents a solid/melt interface \cite{BrosaPlanella2020,Hill1987,McCue2008} or interfacial flow problems where the moving boundary is the interface between two fluids \cite{Crank1987}.  In diffusion-dominated processes, a moving boundary can be driven by swelling of one phase, for example a polymer or a porous material \cite{Mitchell2014,McCue2011,McGuiness2000}. One key motivation from a physics perspective is to use moving boundary problems to describe pattern formation at the interface, such as that which occurs when a crystal forms from a supercooled melt \cite{Gibou2003,Juric1996} or when viscous fingering develops in a Hele-Shaw cell \cite{Howison1986,Langer1989}.  From a mathematical perspective, there is significant interest in studying various asymptotic or singular solutions to these models and their relevance for the underlying physics \cite{Howison1985}.

More recently, researchers have formulated a range of moving boundary problems to model the tissue dynamics and the collective motion of biological cells or bacteria \cite{Franks2003,Giverso2015,Ward2003}.  In this study, we shall be concerned with one of the more simple such moving boundary problems, for which the governing partial differential equation is the Fisher-KPP equation \cite{Du2010}.  Our focus will be on parameter regimes for which the moving boundary retreats towards the population, instead of invading an empty space.  Here the solutions may evolve to travelling waves, undergo extinction, or blow-up in finite time, the latter of which involves the speed of the moving boundary and the slope of the population density both becoming unbounded at the blow-up time.  We shall establish links with ill-posed models for melting/freezing.

In particular, we will study the following dimensionless moving boundary problem
\begin{equation}
\frac{\partial u}{\partial t}=\frac{\partial^2 u}{\partial x^2}+u(1-u),
\quad 0<x<s(t),
\label{eq:StefanFisher1}
\end{equation}
\begin{equation}
\frac{\partial u}{\partial x}=0 \quad\mbox{on}\quad x=0,
\label{eq:StefanFisher2}
\end{equation}
\begin{equation}
u=0, \quad \frac{\mathrm{d}s}{\mathrm{d}t}=-\kappa\frac{\partial u}{\partial x}
\quad\mbox{on}\quad x=s(t),
\label{eq:StefanFisher3}
\end{equation}
\begin{equation}
u(x,0)=F(x), \quad 0<x<s(0).
\label{eq:StefanFisher4}
\end{equation}
Here, $u(x,t)$ represents the density of a certain cell type (such as tumour cells or skin cells), which satisfies the Fisher-KPP equation (\ref{eq:StefanFisher1}) \cite{Fisher1937,Kolmogorov1937}.  In this simple model, cell motility is described by a constant diffusion coefficient, while cell proliferation is governed by the logistic growth term (here the carrying capacity has been scaled to be $u=1$) \cite{Edelstein2005,Murray2002}.  Extensive discussion and analysis of how to apply the Fisher-KPP equation to experiments in cell biology are provided in Refs~\cite{Jin2016,Jin2017,Jin2018,Johnston2015,Maini2004b,Warne2019}, for example.  The domain for the Fisher-KPP equation (\ref{eq:StefanFisher1}) is bounded to the right by a moving boundary $x=s(t)$, at which the cell density vanishes. The second condition in (\ref{eq:StefanFisher3}) is a Stefan-like boundary condition \cite{Hill1987} and for this reason we refer to (\ref{eq:StefanFisher1})-(\ref{eq:StefanFisher4}) as the {\em Fisher-Stefan} model.  The parameter $\kappa$ in (\ref{eq:StefanFisher3}) is a leakage coefficient that drives the speed of the interface in the same way that a Stefan number does in the analogous Stefan problem.  The loss of mass at the moving boundary could be driven by invasion into ($\kappa>0$) or recession from ($\kappa<0$) a harsh environment.  Alternatively, for $\kappa<0$, the cell loss could be associated with a conversion of the main cell population into another cell type that is left behind as the front recedes.  An argument for how $\kappa$ could be determined given experimental data is provided in Ref.~\cite{Elhachem2021}.

The one-phase model (\ref{eq:StefanFisher1})-(\ref{eq:StefanFisher4}) with a single moving boundary sits within a class of Stefan-type moving boundary problems for reaction diffusion models with a linear diffusion term.  Much of this work began with Du \& Lin~\cite{Du2010}, who studied advancing fronts that either go to extinction or evolve to travelling waves.  They refer to these two option as leading to a ``spreading-vanishing dichotomy''.  Rigorous results for this one-phase problem were subsequently established for the long-time asymptotics, more general source terms, and the radially symmetric version in higher dimensions by Du and co-workers~\cite{Du2014a,Du2014b,Du2015}.  Since that time, much attention has been devoted to various generalisations, many of which are reviewed by Du~\cite{Du2020}.  In very recent times, a variety of extensions have been documented, including a rigorous results for different boundary conditions, nonlinear or nonlocal diffusion, generalised reaction terms, with one or two phases, as well as motivation in terms of applications in ecology \cite{Bao2018,Cai2017,Cai2020,Cao2019,Chang2020,Ding2020,Du2018,Endo2020,Jiang2018,Kaneko2018,Lutscher2020,Wang2018,Wang2019,Yuan2021}.  Formal results and numerical simulations have been described by us \cite{Elhachem2019,Elhachem2020,Elhachem2021,McCue2020,Simpson2020} and others \cite{Liu2020}.

Most of the key results for the core Fisher-Stefan model (\ref{eq:StefanFisher1})-(\ref{eq:StefanFisher4}) hold for $\kappa>0$, for which the moving boundary $x=s(t)$ moves in the positive $x$-direction.  For example, with $\kappa>0$ and provided $s(0)>\pi/2$, the time-dependent solution of (\ref{eq:StefanFisher1})-(\ref{eq:StefanFisher4}) will evolve towards a travelling wave profile with a speed $0<c<2$ that depends on $\kappa$, with $c\rightarrow 0^+$ as $\kappa\rightarrow 0^+$ and $c\rightarrow 2^-$ as $\kappa\rightarrow\infty$ \cite{Du2010,Du2020,Elhachem2019}.  This interval of $c$ values for the Fisher-Stefan model is surprising because solutions of the traditional Fisher-KPP equation on $0<x<\infty$ evolve to travelling wave profiles with speeds $c>2$ and the travelling wave profiles for $0<c<2$ (whose trajectories in the phase plane have spiral at the origin) are normally discarded as being unphysical, since they involve negative population densities for sufficiently large values of the independent variable $z=x-ct$.  A stronger result for (\ref{eq:StefanFisher1})-(\ref{eq:StefanFisher4}) is that if $s(t)$ ever becomes greater than $\pi/2-\kappa\int_0^{s(t)}u(x,t)\,\mathrm{d}x$ then the travelling wave is guaranteed \cite{Elhachem2019}.  Alternatively, for $s(0)<\pi/2$ and $u(x,0)=F(x)$ sufficiently small, the solution goes extinct with $u(x,t)\rightarrow 0$, $s(t)\rightarrow s_{\mathrm{e}}$ as $t\rightarrow\infty$, where $s(0)<s_{\mathrm{e}}<\pi/2$ is a constant \cite{Du2010,Du2020,Elhachem2019}. This behaviour is interesting because solutions of the Fisher-KPP equation on $0<x<\infty$ with a no-flux condition at $x=0$ do not ever go extinct.  Together the above results are referred to by Du and coworkers as the spreading-vanishing dichotomy \cite{Du2010,Du2020}.

Much less attention has been devoted to (\ref{eq:StefanFisher1})-(\ref{eq:StefanFisher4}) for $\kappa < 0$ \cite{Elhachem2021,McCue2020}.  Here the moving boundary $x=s(t)$ propagates in the negative $x$-direction.  For $-1<\kappa<0$, we have previously shown that, for sufficiently large $s(0)$, solutions to (\ref{eq:StefanFisher1})-(\ref{eq:StefanFisher4}) very quickly evolve towards receding travelling wave profiles that persist as long as $s(t)\gg 1$ (ultimately the no-flux condition at $x=0$ will disturb this trend)
\cite{Elhachem2021}. These travelling wave solutions move with speed $c\sim \kappa/\sqrt{3}$ as $\kappa\rightarrow 0$ and $c\sim -2^{-1}(1+\kappa)^{-1/2}$ as $\kappa\rightarrow -1^+$; their profiles correspond to portions of trajectories in the phase-plane that emerge from a saddle point (but are not heteroclinic orbits).  An exact solution for $c=-5/\sqrt{6}$ can be written down in terms of the equianharmonic case of the Weierstrass p-function \cite{McCue2020,Ablowitz1979}.  There is a stationary wave solution for $\kappa=0$, which corresponds to a homoclinic orbit in the phase-plane \cite{Elhachem2021}.  Again, all of these travelling wave trajectories for $c\leq 0$ are normally discarded as being unphysical, as their full trajectories in the phase-plane involve negative population densities.  Note the shape of the travelling wave profiles for $-1<\kappa<0$ are qualitatively similar to those for $\kappa>0$, except that they become steeper as $\kappa$ decreases.

In section~\ref{sec:blowup} of the paper, we are mostly concerned with presenting new results for the parameter regime $\kappa<-1$, $s(0)\gg 1$, for which a phase-plane analysis shows there are no travelling wave solutions.  We demonstrate that time-dependent solutions to (\ref{eq:StefanFisher1})-(\ref{eq:StefanFisher4}) can blow up in finite time for $\kappa<-1$, with $\mathrm{d}s/\mathrm{d}t\rightarrow -\infty$ as $t\rightarrow t_\mathrm{c}$, and conjecture that for the borderline case $\kappa=-1$ the blow up is at infinite time.  These results are interesting because this type of blow up has not ever been observed before for models with the Fisher-KPP equation.  For $s(0)=\mathcal{O}(1)$, in section~\ref{sec:stefan} we re-scale the problem in order to compare with a well-studied ill-posed Stefan problem for melting a superheated solid (or freezing a supercooled liquid). In this case, as well as the option of blow-up in finite time, solutions may go to extinction, whereby the population approaches zero in a nonzero domain in the long time limit.  We argue that for $s(0)=\mathcal{O}(1)$, the critical value of $\kappa$ that separates blow-up from extinction depends on the initial population and a parameter that measures cell proliferation (this parameter is equal to $s(0)^2$ in the version of the problem (\ref{eq:StefanFisher1})-(\ref{eq:StefanFisher4})).
Further, we provide estimates for the distance the moving boundary travels before the solution goes to extinction.  As such, we provide new links between the Fisher-Stefan moving boundary model with ill-posed Stefan problems.  Finally, we provide concluding remarks in section~\ref{sec:conclusion}.
Full details of our numerical scheme are provided in Appendix~\ref{sec:Numericalmethods}, while MATLAB software to solve the time-dependent PDE models in this work are available on Github at \href{https://github.com/ProfMJSimpson/Fisher_blowup}{https://github.com/ProfMJSimpson/Fisher\_blowup}.


%


\section{Travelling waves and blow-up for $s(0)\gg 1$}
\label{sec:blowup}

\subsection{Time-dependent solutions}

In this section we shall explore the time-dependent behaviour of solutions to (\ref{eq:StefanFisher1})-(\ref{eq:StefanFisher4}) for cases in which $s(0)\gg 1$.  The motivation here is to observe the solution evolving to a travelling wave profile ($-1<\kappa<0$) or blowing up in finite time ($\kappa<-1$).  Provided that $s(t)\gg 1$, the consequences of imposing the no-flux condition (\ref{eq:StefanFisher2}) are negligible (thus, in effect, we are treating (\ref{eq:StefanFisher1})-(\ref{eq:StefanFisher4}) on $-\infty<x\leq s(t)$).

For that purpose, we use the initial condition
\begin{equation}
F(x)=\alpha (1-\mathrm{H}(s(0))),
\label{eq:IC}
\end{equation}
which is simply a step function, where $\mathrm{H}(x)$ is the Heaviside function and $\alpha>0$ is a parameter that stipulates the initial density.
A crucial step for our numerical scheme is to map $0<x<s(t)$ to a fixed domain $0<y<1$ via $y=x/s(t)$.  The solutions we are interested in turn out to be very steep near $x=s(t)$ ($y=1$) and very flat near $x=0$ ($y=0$), and so it proves useful to discretise the unit interval in $y$ with an uneven mesh such that the grid points are increasing close together as $y$ increases.  For most of the simulations in this paper we use $N=1001$ grid points so that the smallest grid spacing near $y=1$ is $1\times 10^{-6}$.  The details of the numerical scheme are recorded in Appendix~\ref{sec:Numericalmethods}.

In figure~\ref{figure1} we show numerical solutions for three different values of leakage coefficient $\kappa$, all computed with the initial condition (\ref{eq:IC}) with $\alpha=0.5$ and $s(0)=1000$.  In panel (a), the solution for $\kappa=-0.98$ is representative of solutions for $-1<\kappa<0$.  Here, as the interface $x=s(t)$ begins to propagate in the negative $x$-direction, the maximum value of $u$ begins to rise from $u=0.5$ at $t=0$ towards $u=1$, and the solution appears to evolve to a constant shape moving from right to left with constant speed (the solutions are drawn at equally spaced times $t=0$, $1$, $2$, \ldots).  Indeed, as demonstrated in El-Hachem et al.~\cite{Elhachem2021}, this solution is evolving towards a travelling wave profile that we explore further below in section~\ref{sec:twsolns}.  Of course, this convergence to a travelling wave cannot continue indefinitely, as we have imposed a no-flux condition at $x=0$ in our numerical scheme; however, these almost-travelling-wave profiles persist until $s(t)=\mathcal{O}(1)$.  We defer a discussion about the effects of (\ref{eq:StefanFisher2}) until section~\ref{sec:stefan}.  Note that this qualitative behaviour we see in figure~\ref{figure1} for the initial condition with $\alpha=0.5$ is observed for all $0<\alpha<1$.

\begin{figure}
	\centering
	\includegraphics[width=1.0\linewidth]{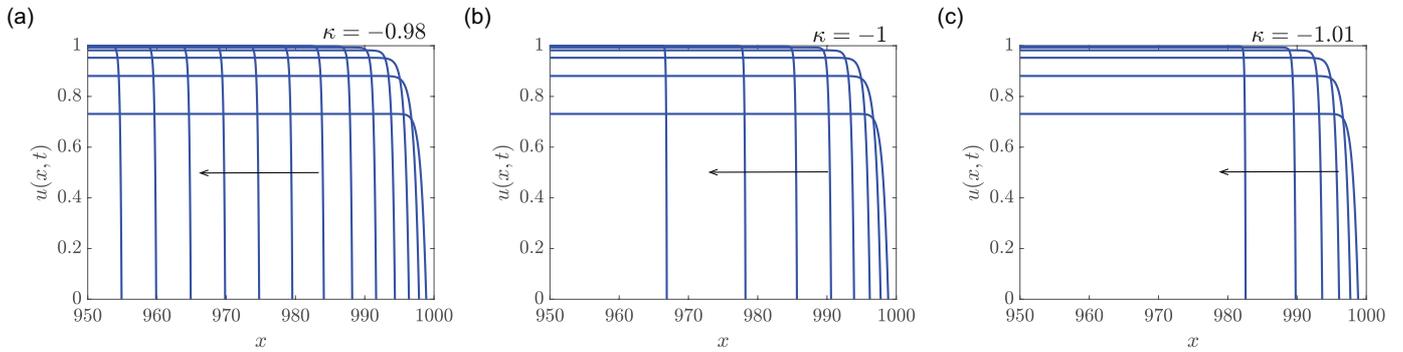}
	\caption{Time-dependent solutions of (\ref{eq:StefanFisher1})-(\ref{eq:StefanFisher4}) for three different values of $\kappa$: (a) $\kappa=-0.98$; (b) $\kappa=-1$; (c) $\kappa=-1.01$.  In each case we have use the initial condition (\ref{eq:IC}) with $\alpha=0.5$ and $s(0)=1000$.  The solutions are drawn at unit time intervals up to (a) $t=13$, (b) $t=8$, (c) $t=6$.  These three examples are representative of solutions that appear to (a) evolve to a travelling wave, (b) blow up in infinite time, (c) undergo finite-time blow-up.}
	\label{figure1}
\end{figure}

Figures \ref{figure1}(b)-(c) contain two more numerical solutions, one for $\kappa=-1$ and the other for $\kappa=-1.01$.  These solutions are for parameter values not considered previously in the literature.  In both cases, the solution does not appear to evolve to a travelling wave.  In (b), with $\kappa=-1$, it is clear the speed of the interface continues to increase as time progresses (again, the solutions are drawn at equally spaced times), while in (c), with $\kappa=-1.01$, this trend is also apparent, although this figure does not show what happens for latter times.  We shall defer further discussion on these solutions for $\kappa\leq -1$ until subsection~\ref{sec:finitetime}.

\subsection{Travelling wave solutions, $-1<\kappa<0$}\label{sec:twsolns}

In this subsection we shall briefly review the key results for the existence of travelling wave solutions for the parameter range $-1<\kappa<0$ \cite{Elhachem2021}.  With this in mind, we write $u=U(z)$, where $z$ is the travelling wave coordinate $z=x-ct$, so that $U$ satisfies
\begin{equation}
\frac{\mathrm{d}^2U}{\mathrm{d}z^2}+c\frac{\mathrm{d}U}{\mathrm{d}z}+U(1-U)=0.
\label{eq:ODEUz}
\end{equation}
Here $c$ is the speed of the travelling wave solution.

For the physically relevant boundary conditions we have in mind here, namely
\begin{equation}
U\rightarrow 1^-, \quad \frac{\mathrm{d}U}{\mathrm{d}z}\rightarrow 0^-,
\quad\mbox{as}\quad z\rightarrow -\infty,
\label{eq:Urightarrow1}
\end{equation}
most studies of (\ref{eq:ODEUz}) are for $c>0$.   In that case, traditional travelling wave solutions to the Fisher-KPP equation, for $c\geq 2$, correspond to heteroclinic orbits that join the saddle point representing (\ref{eq:Urightarrow1}) with the origin $U=0$, $\mathrm{d}U/\mathrm{d}z=0$.  As mentioned in the Introduction, for the Fisher-Stefan problem (\ref{eq:StefanFisher1})-(\ref{eq:StefanFisher4}) with $\kappa>0$ and sufficiently large initial values of $s(0)$, solutions evolve to travelling wave solutions for $0<c<2$; these correspond to the portion of the heteroclinic orbit joining the saddle and the origin that has $U\geq 0$, $\mathrm{d}U/\mathrm{d}z\geq 0$.

Turning now to the regime we are interested in here, namely $c<0$.  In this case we let $V=-\mathrm{d}U/\mathrm{d}z$, so that (\ref{eq:ODEUz}) is rewritten as the coupled system of first order equations
\begin{align}
\frac{\text{d}U}{\text{d} z} & = -V, \label{eq:ODEdU}\\
\frac{\text{d}V}{\text{d} z} &= -cV + U(1-U), \label{eq:ODEdV}
\end{align}
The phase-plane $(U,V)$ is simply a reflection of the traditional representation about the $U$-axis.  As such, the qualitative behaviour of trajectories carries over from before.  That is, for $-2<c<0$ the equilibrium point at the origin is a spiral, just like that for $0<c<2$ except for the reflection, while the origin for $c\leq -2$ is a node as with $c\geq 2$, again with reflection.

\begin{figure}
	\centering
	\includegraphics[width=.8\linewidth]{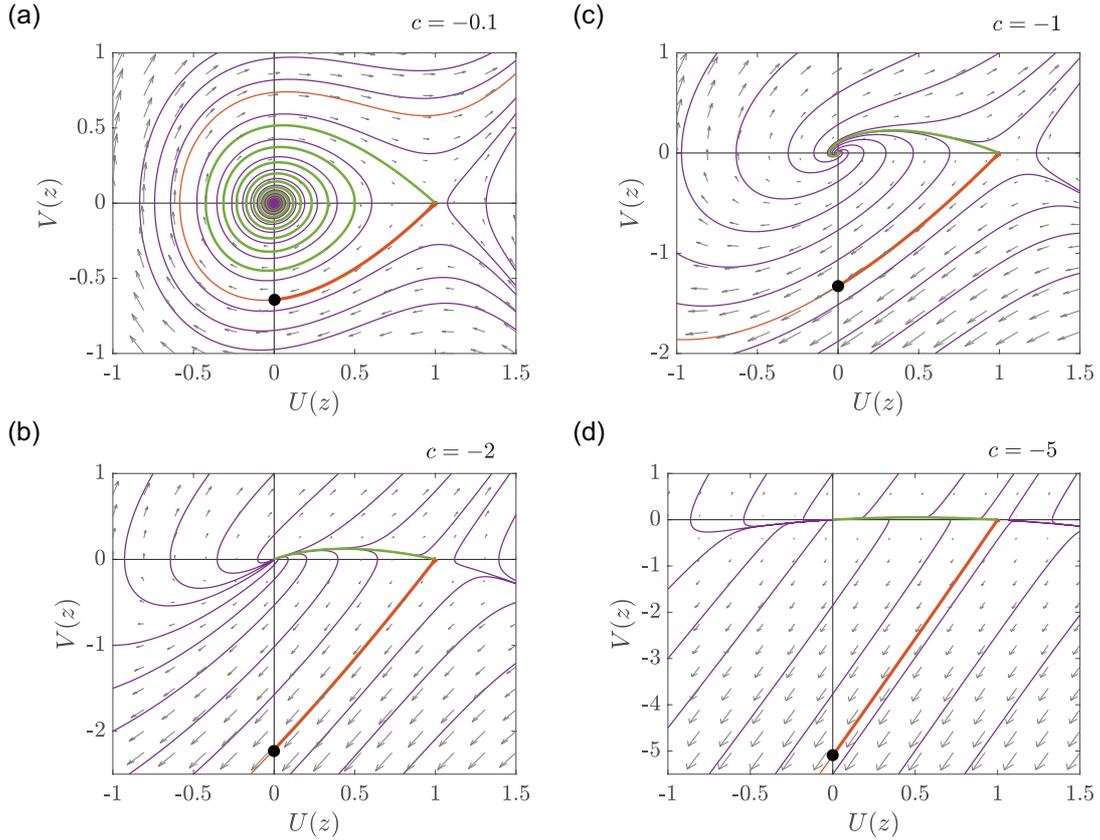}
	\caption{Representative phase-planes for (a) $c=-0.1$, (b) $-1$, (c) $-2$ and (d) $-5$, drawn by numerically integrating (\ref{eq:ODEdU})-(\ref{eq:ODEdV}).  In each case, the thick trajectory in red represents the solution to (\ref{eq:ODEUz}) that satisfies the far-field condition (\ref{eq:Urightarrow1}).  These trajectories are truncated at the $V$-axis at some value $V=V^*$ (indicated by the solid black disc) which allows for the numerical prediction of $\kappa$ via $c=-\kappa V^*$, thereby linking these travelling wave profiles with the intermediate-time solutions to the moving boundary problem (\ref{eq:StefanFisher1})-(\ref{eq:StefanFisher4}).  Note the green thick trajectories are the heteroclinic orbits that connect the saddle point $(U,V)=(1,0)$ with the fixed point at the origin.}
	\label{figure2}
\end{figure}

We demonstrate these ideas with four examples in figure~\ref{figure2}, drawn for the different values of $c=-0.1$, $-1$, $-2$ and $-5$.  These results are found by solving (\ref{eq:ODEdU})--(\ref{eq:ODEdV}) numerically using Heun's method with a constant step size $\textrm{d}z$.  In each case shown in this figure, we have drawn (thick, red) the portions of the phase-plane that correspond to the travelling wave solutions we are concerned with.  These all begin at $(U,V)=(0,0)$ and traverse the fourth quadrant until they intersect the $V$-axis at some special value, $V=V^*$ say.  For the four examples presented in figure~\ref{figure2}, the numerical estimates of $V^*$ are (a) $c=-0.1$, $V^*=-0.643$, (b) $c=-1$, $V^*=-1.32$, (c) $c=-2$, $V^*=-2.21$ and (d) $c=-5$, $V^*=-5.10$.  Given the Stefan condition (\ref{eq:StefanFisher3}), we know that $c=-\kappa V^*$, so these numerical estimates in the phase plane correspond to (a) $\kappa=-0.156$, (b), $\kappa=-0.753$, (c) $\kappa=-0.904$ and (d) $\kappa=-0.981$.  For each of these $\kappa$ values, our numerical simulations show that time-dependent solutions of the full Fisher-Stefan model (\ref{eq:StefanFisher1})-(\ref{eq:StefanFisher4}) with $s(0)\gg 1$ evolve towards the travelling wave profiles with the predicted values of $c$.

For completeness, we also show in figure~\ref{figure2} a number of trajectories (thin, purple) that do not begin at the saddle point at $(U,V)=(1,0)$ and so do not satisfy (\ref{eq:Urightarrow1}).  None of these are physically-relevant travelling wave solutions.  We have also shown the trajectories that represent either invading travelling wave solutions (thick, green).  For (a)-(b), for which the values of $c$ lie in the interval $-2<c<0$, these trajectories are normally discarded since the spirals involve negative portions of the population; however, for the Fisher-Stefan model (\ref{eq:StefanFisher1})-(\ref{eq:StefanFisher4}) with $\kappa>0$, these are truncated at when they first meet the $V$-axis.  For (c)-(d), for which $c\leq -2$, the invading travelling wave trajectories (thick, green) represent traditional travelling wave solutions to the Fisher-KPP equation.

It is easy to see that the solutions of (\ref{eq:ODEUz}) for $c<0$ that satisfy (\ref{eq:Urightarrow1}) (corresponding to the thick red curves in figure~\ref{figure2}) exist for $-\infty<c<0$.  In the limit $c\rightarrow -\infty$, we see from (\ref{eq:ODEdU})-(\ref{eq:ODEdV}) that $\mathrm{d}V/\mathrm{d}U=-c+U(1-U)/V \sim -c$ which suggests the relevant trajectories in the phase-plane approach the straight lines $V=c(U-1)$.  Given $V^*=\left. V\right|_{U=0}$, we have $\kappa\rightarrow -1^+$ as $c\rightarrow -\infty$.  In El-Hachem et al.~\cite{Elhachem2021} we derive the approximation
\begin{equation}\label{eq:Asymptotic}
U\sim 1-\mathrm{e}^{-cz}+\frac{\mathrm{e}^{-cz}\left(2cz-1+\mathrm{e}^{-cz}\right)}{2c^2}
\quad\mbox{as}\quad
c\rightarrow -\infty,
\end{equation}
which leads to the limiting behaviour $\kappa\sim -1+(2c^2)^{-1}$ as $c\rightarrow -\infty$
or, in other words, $c\sim -2^{-1}(1+\kappa)^{-1/2}$ as $\kappa\rightarrow -1^+$.

\subsection{Finite-time blow up, $\kappa<-1$}\label{sec:finitetime}

The previous subsection summarises recently-recorded results for values of the leakage coefficient in the interval $-1<\kappa<0$, for which time-dependent solutions of (\ref{eq:StefanFisher1})-(\ref{eq:StefanFisher4}) evolve to travelling waves.  Since the phase-plane analysis in that subsection does not apply for $\kappa\leq -1$, we are now motivated to explore the time-dependent behaviour of solutions in this parameter regime.

Returning to figure~\ref{figure1}, the solutions in (b) and (c) do not appear to approach a travelling wave like the solution in (a).  For example, the moving boundary in (b) ($\kappa=-1$) appears to move a longer distance at each time-interval, remembering that each solution is drawn at unit intervals $t=0$, $1$, $2$, \ldots.  The largest time shown in figure~\ref{figure1}(b) is $t=7$ and indeed the location of the moving boundary for $t=8$ is too far to the left to fit on the figure.  Similarly, for the solution in panel (c) ($\kappa=-1.01$), the largest time shown is $t=5$ since the solution for $t=6$ does not fit on this scale.

\begin{figure}
	\centering
	\includegraphics[width=1\linewidth]{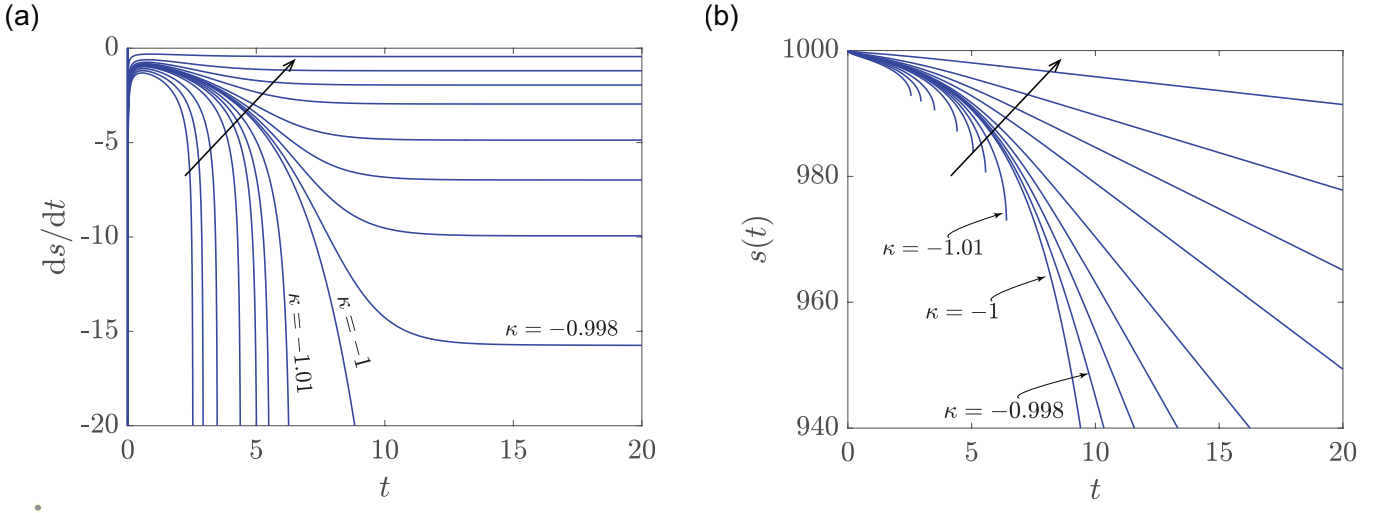}
	\caption{Numerical results for (\ref{eq:StefanFisher1})-(\ref{eq:StefanFisher4}) showing the dependence of (a) the speed of the moving boundary $\mathrm{d}s/\mathrm{d}t$ and (b) location of the moving boundary $s(t)$ versus $t$.  The arrows indicate increasing values of $\kappa$.  From bottom left to top right, the solutions are drawn for $\kappa=-1.2$, $-1.15$, $-1.1$, $-1.05$, $-1.03$, $-1.02$, $-1.01$, $-1$, $-0.998$, $-0.995$, $-0.99$, $-0.98$, $-0.95$, $-0.9$, $-0.8$ and $-0.5$.}
	\label{figure3}
\end{figure}

To further explore the ultimate behaviour of these solutions, we present in figure~\ref{figure3}(a) plots of the speed of the moving boundary $x=s(t)$ versus time for various values of $\kappa$.  For the eight examples with $-1<\kappa<0$, we clearly see the speed levelling off to a constant value as time increases, supporting our observations that solutions in this parameter regime evolve to travelling waves.  On the other hand, for the six examples with $\kappa<-1$, the speed appears to approach a vertical asymptote, suggesting there is a form of finite-time blow-up with $\mathrm{d}s/\mathrm{d}t\rightarrow -\infty$ as $t\rightarrow t_{\mathrm{c}}^-$.  For instance, for $\kappa=-1.01$, corresponding to the solution in figure~\ref{figure1}(c), the blow-up time appears to be roughly $t_{\mathrm{c}}\approx 6.4$.

A slightly different perspective is provided by plotting the moving boundary location $s(t)$ versus time $t$, as in figure~\ref{figure3}(b).  Here we see solutions for $-1<\kappa<0$ evolve to straight lines with nonzero slopes that correspond to the travelling wave speeds, while solutions for $\kappa<-1$ appear to truncate at finite values of $s$ at the finite blow-up time, which suggests that $s\rightarrow s_{\mathrm{c}}^+$ as $t\rightarrow t_{\mathrm{c}}^-$.  For the particular example $\kappa=-1.01$, the numerical estimate for the end point is $s_{\mathrm{c}}\approx 975$.

\begin{figure}
	\centering
	\includegraphics[width=1\linewidth]{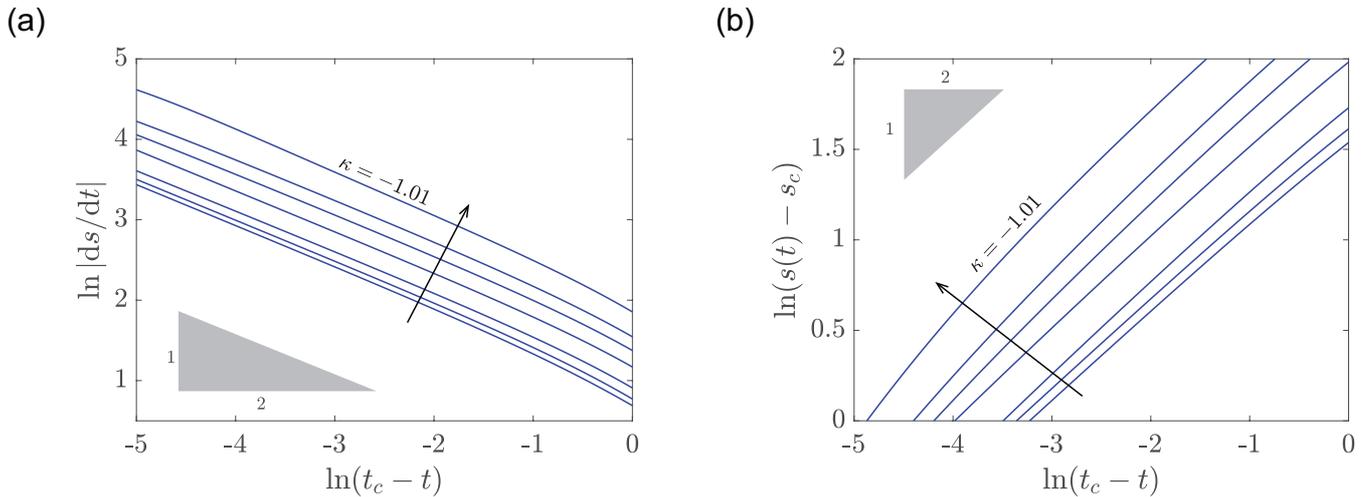}
	\caption{Log-log plots for (a) $\mathrm{d}s/\mathrm{d}t$ versus $t_{\mathrm{c}}-t$ and (b) $s(t)-s_{\mathrm{c}}$ versus $t_{\mathrm{c}}-t$, drawn for $\kappa=-1.2$, $-1.15$, $-1.1$, $-1.05$, $-1.02$, $-1.01$.  Also included is a line segment with slope $-1/2$.}
	\label{figure4}
\end{figure}

It is useful to explore the scalings in the blow-up limit $t\rightarrow t_{\mathrm{c}}^-$.  We plot in figure~\ref{figure4} the dependence of $\ln(-\mathrm{d}s/\mathrm{d}t)$ and $\ln(s(t)-s_{\mathrm{c}})$ versus $\ln(t_{\mathrm{c}}-t)$.  On this scale, it appears that in both (a) and (b) the curves are roughly linear with slope $-1/2$, suggesting a power-law relationship in the limit.  We return to this issue in Section~\ref{sec:stefan}.

\begin{figure}
	\centering
	\includegraphics[width=0.5\linewidth]{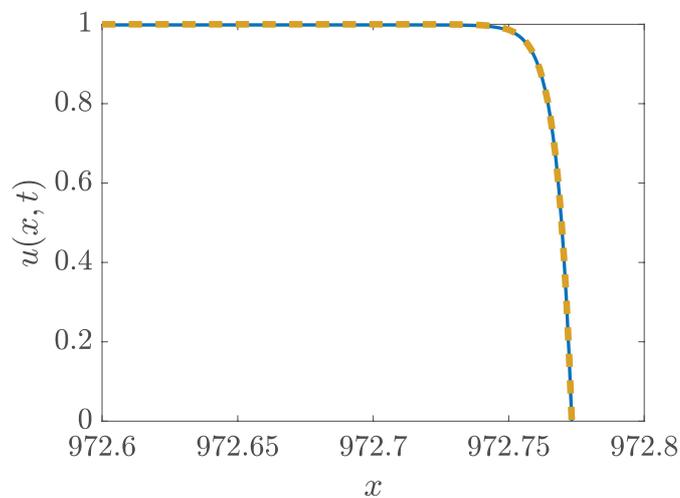}
    \caption{Numerical solution very close to blow-up.  The solid blue curve is the numerical solution of (\ref{eq:StefanFisher1})-(\ref{eq:StefanFisher4}) for $\kappa=-1.01$, $t=6.41$, for which $s\approx 972.78$ and $\dot{s}\approx -183.78$.  The dashed orange curve is the estimate (\ref{eq:crude}) computed using these values.}
    \label{newfigure5}
\end{figure}

The shape of the solution at blow-up appears to be qualitatively similar to that for times leading up to blow-up.  For example, for the simulation shown in  figure~\ref{figure1}(c), where we have $\kappa=-1.01$, the numerical estimate of the blow-up profile is included in figure~\ref{newfigure5} as a solid blue curve.  This numerical solution is computed very slightly before the estimated blow-up time, which is $t_{\mathrm{c}}\approx 6.4$.  While it is difficult to appreciate the difference between this solution and those for previous times, the key point is that this form of blow-up is characterised by the slope of the solution becoming infinite at $x=s(t)$ where $u(s(t),t)=0$.  A simple estimate of the solution near blow-up for $\kappa\sim -1^-$ can be written down by looking for a solution $u\sim w(\zeta)$, where $\zeta=\dot{s}(x-s(t))$, assuming $|\dot{s}^2|\gg|\ddot{s}|\gg 1$ (which is true if $s(t)-s_{\mathrm{c}}\sim (t_{\mathrm{c}}-t)^{1/2}$ or similar, for example).  Then, to leading order, (\ref{eq:StefanFisher1}) gives $-w'\sim w''$, whose solution is $w= A+B\,\mathrm{e}^{-\zeta}$, where $A$ and $B$ are constants.  By fixing $w=0$ at $\zeta=0$ and matching back to $w=1$ as $\zeta\rightarrow-\infty$, we find
\begin{equation}
u\sim 1-\mathrm{e}^{-\dot{s}(x-s(t))}
\quad \mbox{as}\quad t\rightarrow t_{\mathrm{c}}^-, \quad \kappa\sim -1^-.
\label{eq:crude}
\end{equation}
We plot this simple profile in figure~\ref{newfigure5} using the numerical values of $\dot{s}$ and $s$ for this particular value of $t$.  The agreement is quite good.

Having established strong numerical evidence for travelling wave outcomes for $-1<\kappa<0$ and finite-time blow-up for $\kappa<-1$, it remains to consider the borderline case $\kappa=-1$.  In this instance, we conjecture that the solution undergoes infinite-time blow-up.  As mentioned above, in figure~\ref{figure1}(b) we see the solution front appearing to increase its speed with time, and that observation is confirmed in figure~\ref{figure3}, where both $|\mathrm{d}s/\mathrm{d}t|$ and $|s(t)|$ appear to continue to increase with time.  These numerical results support infinite-time blow-up.  Furthermore, from our phase-plane analysis of travelling wave solutions in Section~\ref{sec:twsolns} and El-Hachem et al.~\cite{Elhachem2021}, we find the travelling wave speeds has the limiting behaviour $c\sim -2^{-1}(1+\kappa)^{-1/2}$ as $\kappa\rightarrow -1^+$, which provides further analytical evidence for our conjecture.

\section{Blow-up versus extinction for $s(0)=\mathcal{O}(1)$}
\label{sec:stefan}

\subsection{Preamble}\label{sec:3preamble}

We turn now to studying (\ref{eq:StefanFisher1})-(\ref{eq:StefanFisher4}) for cases in which $s(0)=\mathcal{O}(1)$.  Here the no-flux condition (\ref{eq:StefanFisher2}) is important and cannot be ignored.

It proves convenient to rescale our problem by introducing the new variables $\tilde{x}=x/s(0)$, $\tilde{t}=t/s(0)^2$, $\tilde{s}(\tilde{t})=s(t)/s(0)$, so that (\ref{eq:StefanFisher1})-(\ref{eq:StefanFisher4}) becomes
\begin{equation}
\frac{\partial \tilde{u}}{\partial \tilde{t}}=\frac{\partial^2 \tilde{u}}{\partial \tilde{x}^2}+\lambda\tilde{u}(1-\tilde{u}),
\quad 0<\tilde{x}<\tilde{s}(\tilde{t}),
\label{eq:StefanFisherscaled1}
\end{equation}
\begin{equation}
\frac{\partial \tilde{u}}{\partial \tilde{x}}=0 \quad\mbox{on}\quad \tilde{x}=0,
\label{eq:StefanFisherscaled2}
\end{equation}
\begin{equation}
\tilde{u}=0, \quad \frac{\mathrm{d}\tilde{s}}{\mathrm{d}\tilde{t}}=-\kappa\frac{\partial \tilde{u}}{\partial \tilde{x}}
\quad\mbox{on}\quad \tilde{x}=\tilde{s}(\tilde{t}),
\label{eq:StefanFisherscaled3}
\end{equation}
\begin{equation}
\tilde{u}(\tilde{x},0)=\tilde{F}(\tilde{x}), \quad 0<\tilde{x}<1,
\label{eq:StefanFisherscaled4}
\end{equation}
where the two parameters in the problem are now $\kappa$ (as before) and
\begin{equation}
\lambda=s(0)^2.
\label{eq:lambda}
\end{equation}
Note that the spatial domain is now a subset of $[0,1]$ so that for $\lambda=\mathcal{O}(1)$ the no-flux condition (\ref{eq:StefanFisherscaled2}) plays an important role.  In particular, in this parameter regime, solutions do not evolve to travelling waves, but instead appear to either undergo finite-time blow-up or go to extinction, as we shall demonstrate.

\subsection{Model for melting a superheated solid}\label{sec:superheated}

The reformulation in section~\ref{sec:3preamble} has the advantage of setting up a direct analogy with a well-studied problem for melting a superheated solid (or freezing a supercooled liquid), as we now summarise.  Consider (\ref{eq:StefanFisherscaled1})-(\ref{eq:StefanFisherscaled4}) with $\lambda=0$.  In this case, the moving boundary problem has the following interpretation \cite{Howison1985,Back2010,Fasano1977,Fasano1981,Herrero1996,King2005,Lacey1985,Sherman1970}.  First, $\tilde{u}(\tilde{x},\tilde{t})$ is the (dimensionless) temperature of a superheated solid on the domain $0<\tilde{x}<\tilde{s}(\tilde{t})$, where here the melting temperature is scaled to be $\tilde{u}=0$.  Here $\tilde{x}=\tilde{s}(\tilde{t})$ is the interface between the solid phase $0<\tilde{x}<\tilde{s}(\tilde{t})$ and the liquid phase $\tilde{x}>\tilde{s}(\tilde{t})$, the latter being assumed to remain at the melting temperature for all time.  The first condition in (\ref{eq:StefanFisherscaled3}) enforces continuity of temperature across the domains, while the second is a Stefan condition that balances heat conducting in to the interface with energy required to melt the solid.  Remembering that $\kappa<0$, this balance is characterised by the Stefan number $\beta=-1/\kappa$, which is a ratio of the latent heat per unit mass at the equilibrium temperature with a product of the specific heat of the solid and a representative temperature scale.

The one-phase Stefan problem (\ref{eq:StefanFisherscaled1})-(\ref{eq:StefanFisherscaled4}) with $\lambda=0$ is nonclassical because the initial temperature $\tilde{F}(\tilde{x})$ is positive, which means the solid is superheated.  Here the solid-melt interface $\tilde{x}=\tilde{s}(\tilde{t})$ is driven to move in the negative $\tilde{x}$-direction as the solid melts.  Such a scenario is particularly unstable and so this simplified model, which is ill-posed, is highly idealised.  In the classical case, for which the initial temperature is negative, the interface would move in the positive $\tilde{x}$-direction.  Note that by replacing $\tilde{u}$ with $-\tilde{u}$, the model (\ref{eq:StefanFisherscaled1})-(\ref{eq:StefanFisherscaled4}) with $\lambda=0$ is for freezing a supercooled liquid, such as that which occurs with diffusion-limited dendritic solidification.  Therefore, the behaviour of the interface $\tilde{x}=\tilde{s}(\tilde{t})$ that we describe for the superheated solid also applies for the equivalent problem for a supercooled liquid.

It is known that there are three qualitatively different outcomes for (\ref{eq:StefanFisherscaled1})-(\ref{eq:StefanFisherscaled4}) with $\lambda=0$ \cite{Howison1985,Back2010,Fasano1977,Fasano1981,Lacey1985,Sherman1970}.  The first is incomplete melting of the solid, whereby the solution exists for all time with $\tilde{s}\rightarrow \tilde{s}_{\mathrm{e}}^+$ and $\tilde{u}(\tilde{x},\tilde{t})\rightarrow 0^+$ as $\tilde{t}\rightarrow\infty$.  The second is complete melting of the solid, where $\tilde{s}(\tilde{t}_\mathrm{e})=0$ for a finite time $\tilde{t}_\mathrm{e}$.  Finally, the third possible outcome is finite-time blow-up of the type we have studied in subsection~\ref{sec:finitetime}, whereby the solution exists up to some time $\tilde{t}_\mathrm{c}$, with $\tilde{s}\rightarrow\tilde{s}_{\mathrm{c}}^+$ and $\mathrm{d}\tilde{s}/\mathrm{d}\tilde{t}\rightarrow -\infty$ as $\tilde{t}\rightarrow\tilde{t}_{\mathrm{c}}^-$.  To help classify possible solutions, the superheating parameter
\begin{equation}
Q = \int \limits_{0}^{1} \! \tilde{F}(\tilde{x})\, \mathrm{d}\tilde{x}+\kappa^{-1}
\label{eq:intro:01Q}
\end{equation}
is useful.  Here $Q$ is the sum of the amount of heat that is required to be removed from the solid in order to reduce the temperature from $\tilde{u}=\tilde{F}(\tilde{x})$ to $\tilde{u}=0$ and the latent heat which must be absorbed at the interface in order to melt the solid (remembering that no heat can enter the system from the left-hand boundary $\tilde{x}=0$).  A simple result is that the third qualitative outcome (finite-time blow-up) will always occur if $Q>0$, as the initial heat energy in the solid is greater than what is required to melt it, leaving a surplus of energy which cannot conduct away or leave the system.  Another is if $\tilde{F}(\tilde{x})$ is smooth and monotonic with $\tilde{F}(1)=0$, then $Q<0$ leads to the first qualitative outcome, as there is not enough heat initially in the solid to convert to the latent heat required to melt the entire solid, while $Q=0$ is the borderline case in which there is precisely the right amount of initial heat to melt the entire solid and the region melts completely in finite time.

\begin{figure}
	\centering
	\includegraphics[width=1\linewidth]{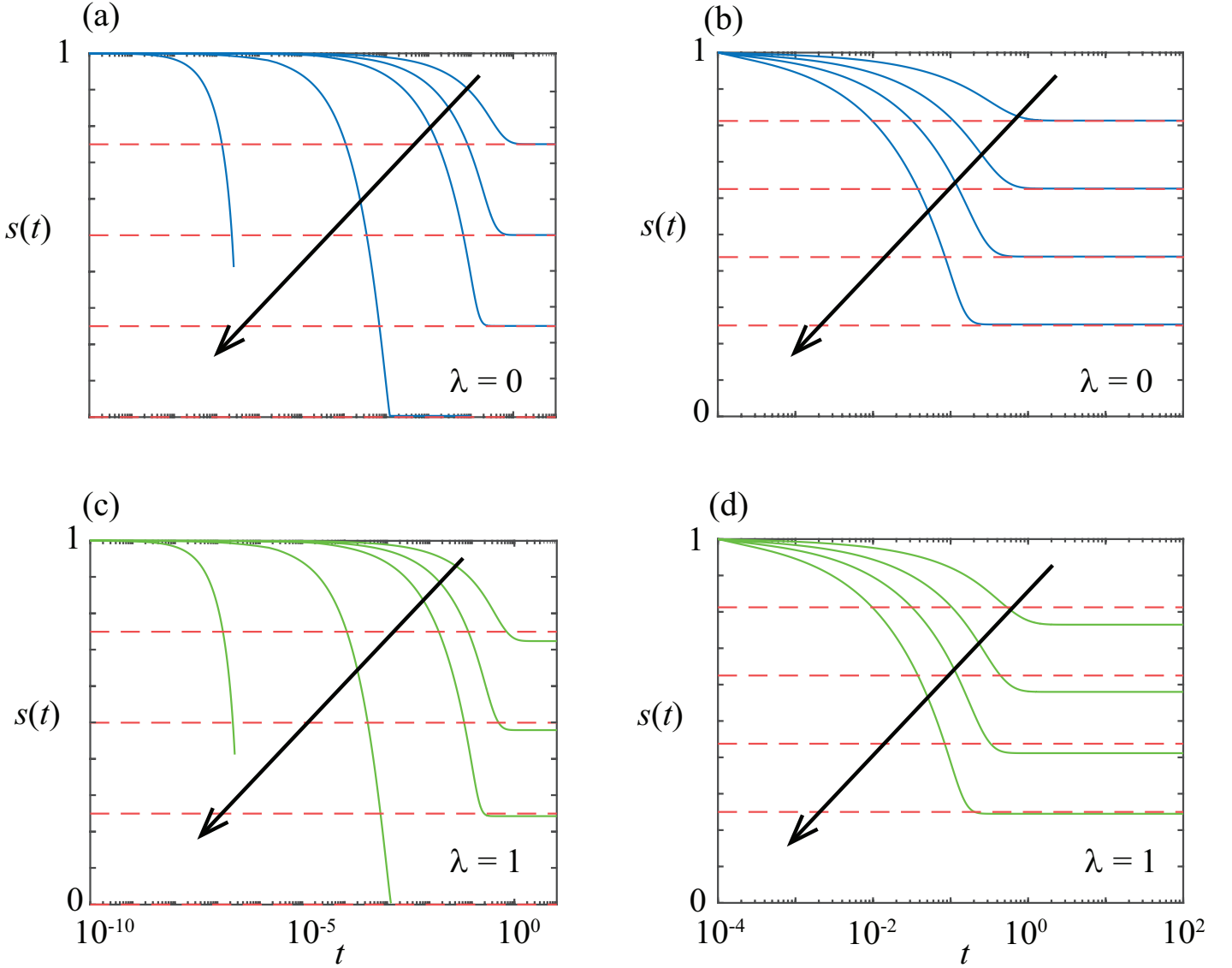}
	\caption{Numerical solutions of (\ref{eq:StefanFisherscaled1})-(\ref{eq:StefanFisherscaled4}) presented in terms of $\tilde{s}(\tilde{t})$ with for $\tilde{s}(0)=1$. (a) $M(0)=1$ with $\kappa=-0.25$ ($Q=-3$), $-0.5$ ($Q=-1$), $-0.75$ ($Q=-0.333$), $-1$ ($Q=0$) and $-1.25$ ($Q=0.2$), with the arrow showing the direction of decreasing $\kappa$. (b) $\kappa = -0.75$ with $M(0)=1$ ($Q=-0.333$), 0.75 ($Q=-0.583$), 0.5 ($Q=-0.833$) and 0.25 ($Q=-1.083$), with the arrow showing the direction of increasing $M(0)$.  Results in (c)--(d) are the same as in (a)--(b), respectively, except that $\lambda=1$.  In each subfigure we show a red dashed line indicating $\tilde{s}_{\mathrm{e}}=\kappa Q$, which is the predicted final position of the moving boundary when $\lambda=0$ (except for the cases $M(0)=1$, $\kappa=-1.25$, for which the prediction $\tilde{s}_{\mathrm{e}}<0$).}
	\label{figure6}
\end{figure}

In order to explore these results further, we integrate (\ref{eq:StefanFisherscaled1}) with $\lambda=0$ and combine with (\ref{eq:StefanFisherscaled2})-(\ref{eq:StefanFisherscaled4}) to give
\begin{equation}
\tilde{s}(\tilde{t})=-\kappa(M(\tilde{t})-Q),
\label{eq:prection1}
\end{equation}
where $M(\tilde{t})$ is the amount of heat in the solid at a given time,
$$
M(\tilde{t})=\int_0^{\tilde{s}(\tilde{t})}\,\tilde{u}(\tilde{x},\tilde{t}) \, \mathrm{d}\tilde{x}.
$$
Then for $Q<0$ with incomplete melting, whereby $\tilde{s}\rightarrow \tilde{s}_{\mathrm{e}}^+$ and $\tilde{u}(\tilde{x},\tilde{t})\rightarrow 0^+$ as $\tilde{t}\rightarrow\infty$, we can predict that $\tilde{s}_{\mathrm{e}}=\kappa Q$ or, alternatively
\begin{equation}
\tilde{s}_{\mathrm{e}}=1+\kappa M(0),
\quad\mbox{for}\quad
Q<0\quad \mbox{(incompete melting)}
\label{complete}
\end{equation}
\cite{Howison1985,Back2010,Fasano1977,Fasano1981,Herrero1996,King2005,Lacey1985,Sherman1970}.
We demonstrate this behaviour in figure~\ref{figure6} where we show time-dependent solutions of (\ref{eq:StefanFisherscaled1})-(\ref{eq:StefanFisherscaled4}) with $\lambda=0$ for an initial condition $\tilde{u}(\tilde{x},0)=M(0)$, $0\leq\tilde{x}\leq 1$.  In figure~\ref{figure6}(a), we have set $M(0)=1$ and used the leakage coefficient values $\kappa=-0.25$, $-0.5$, $-0.75$, $-1$ and $-1.25$.  Here the predicted final positions of the moving boundary are $\tilde{s}_{\mathrm{e}}=0.75$, $0.5$, $0.25$, $0$ and $-0.25$, respectively.  For the first four of these cases, the numerical simulations agree well with the analytical prediction.  For the last case, $\kappa=-1.25$, the prediction $\tilde{s}_{\mathrm{e}}=-0.25$ is unphysical; however, here $Q>0$ and so instead the solution undergoes finite-time blow-up, as described above.  Similarly, in figure~\ref{figure6}(b), we show results for $\kappa=-0.75$ with $M(0)=1$, 0.75, 0.5 and 0.25.  In this case the predicted values of $\tilde{s}_{\mathrm{e}}$ are $0.25$, $0.4375$, $0.625$ and $0.8125$, respectively.  Again, the numerical simulations agree well with the analytical results.

\subsection{Time-dependent results for $\lambda>0$}\label{sec:skellam}

The case in which $\lambda > 0$ is fundamentally different because the logistic growth term in (\ref{eq:StefanFisherscaled1}) acts to increase the total population, making the full problem (\ref{eq:StefanFisherscaled1})-(\ref{eq:StefanFisherscaled4}) with $\lambda>0$ more complicated than that with $\lambda=0$, especially in terms of analytical predictions.  Results in figure \ref{figure6}(c)--(d) are presented in the same format as those in figure \ref{figure6}(a)--(b) except that we set $\lambda = 1$.  Here we see that the formula for $\tilde{s}(\tilde{t})$ given by (\ref{eq:prection1}) (which is for $\lambda=0$) underestimates the distance interface moves.  This observation agrees with our intuition that the source term in (\ref{eq:StefanFisherscaled1}) increases the total population, so the front moves further towards the origin before extinction.

As a first step towards understanding the extinction behaviour with $\lambda > 0$ we approximate (\ref{eq:StefanFisherscaled1})-(\ref{eq:StefanFisherscaled4}) by neglecting the quadratic term, so that we are working with an analogue of Skellam's equation~\cite{Skellam1951}.  In this case, assuming that the population goes to extinction (that is, $\tilde{s}\rightarrow \tilde{s}_{\mathrm{e}}^+$ and $\tilde{u}(\tilde{x},\tilde{t})\rightarrow 0^+$ as $\tilde{t}\rightarrow\infty$), then we obtain
\begin{equation}
\tilde{s}_{\mathrm{e}}=\kappa \left( Q + \lambda \int_{0}^{\infty}M(t) \, \textrm{d}t \right)
\quad\mbox{for extinction}.
\label{eq:prection2}
\end{equation}
Since the integral is positive and $\kappa < 0$, this expression confirms that $\tilde{s}_{\mathrm{e}}$ decreases with $\lambda$, consistent with our numerical explorations in figure \ref{figure6}, where we compared results for $\lambda=0$ with those for $\lambda=1$.  Interestingly, for the case $M(0)=1$, $\kappa=-1$ in figure~\ref{figure6}(a) ($\lambda=0$), we have the exact prediction $\tilde{s}_{\mathrm{e}}=0$, while for $M(0)=1$, $\kappa=-1$ in figure~\ref{figure6}(c) ($\lambda=1$), our approximate result (\ref{eq:prection2}) has $\tilde{s}_{\mathrm{e}}<0$, which is unphysical and implies the solution blows up in finite time.  This result demonstrates how increasing $\lambda$ results in a transition from extinction to blow-up.

\subsection{Blow-up revisited}\label{sec:blowuprevisited}

As discussed above, for the Stefan problem (\ref{eq:StefanFisherscaled1})-(\ref{eq:StefanFisherscaled4}) with $\lambda=0$, the exact result (\ref{complete}) holds provided $Q<0$ and the solutions goes to extinction ($\tilde{u}\rightarrow 0$ as $\tilde{t}\rightarrow\infty$ for all $0<\tilde{x}<\tilde{s}_{\mathrm{e}}$).  Similarly, the borderline case occurs when $Q=0$ and $\tilde{s}_{\mathrm{e}}=0$.  On the other hand, for $Q>0$, finite-time blow-up must occur.  For the Fisher-Stefan model (\ref{eq:StefanFisherscaled1})-(\ref{eq:StefanFisherscaled4}) with $\lambda>0$, we use our observations from subsection~\ref{sec:skellam} to
conjecture that finite-time blow-up occurs for $Q>Q_{\mathrm{crit}}$, where $Q_{\mathrm{crit}}<0$ depends on $\lambda$, with $Q_{\mathrm{crit}}\rightarrow 0$ as $\lambda\rightarrow 0^-$.

For example, the solution in figure~\ref{figure6}(c) computed for $M(0)=1$ and $\kappa=-1$ has $Q=0$.  Since $\lambda>0$ here, our conjecture here is that blow-up will occur if $Q>Q_{\mathrm{crit}}$ where $Q_{\mathrm{crit}}<0$, including the value $Q=0$.  Therefore this numerical solution is an example of this type of behaviour.

We have already discussed finite-time blow-up in subsection~\ref{sec:finitetime}.  It is worth revisiting the near-blow-up behaviour in some detail.  For the Stefan problem (\ref{eq:StefanFisherscaled1})-(\ref{eq:StefanFisherscaled4}) with $\lambda=0$, it is known that
\begin{equation}
\tilde{s}(\tilde{t})-\tilde{s}_{\mathrm{c}} \sim 2\left(\tilde{t}_{\mathrm{c}}-\tilde{t}\right)^{1/2}\mathrm{ln}^{1/2}\left( -\mathrm{ln}\left( \tilde{t}_{\mathrm{c}}-\tilde{t} \right) \right)\quad \mathrm{as} \quad \tilde{t}\rightarrow \tilde{t}_{\mathrm{c}}^-,
\label{eq:blowup}
\end{equation}
which is very close to a $1/2$ power-law dependence but moderated by a very weak logarithm~\cite{Herrero1996}. This complicated scaling can be derived using formal asymptotics via a Biaocchi transform~\cite{King2005}.  While this machinery is not appropriate for (\ref{eq:StefanFisherscaled1})-(\ref{eq:StefanFisherscaled4}) with $\lambda>0$, it is reasonable to postulate that (\ref{eq:blowup}) holds for $\lambda>0$ as well as $\lambda=0$, since $\partial^2\tilde{u}/\partial\tilde{x}^2\gg \tilde{u}(1-\tilde{u})$ in the neighbourhood of $\tilde{x}=\tilde{s}(\tilde{t})$ close to the blow-up time.

To support our conjecture that (\ref{eq:blowup}) is appropriate for (\ref{eq:StefanFisherscaled1})-(\ref{eq:StefanFisherscaled4}) with $\lambda>0$, a brief summary of the formal analysis without the Baiocchi transform is as follows.  Clearly, $\tilde{u}^2\ll \tilde{u}$ in the neighbourhood of $\tilde{x}=\tilde{s}(\tilde{t})$, thus the important near blow-up behaviour is driven by the solution to
\begin{equation}
\frac{\partial \tilde{u}}{\partial \tilde{t}}=\frac{\partial^2 \tilde{u}}{\partial \tilde{x}^2}+\lambda\tilde{u},
\quad \tilde{s}(\tilde{t})-\tilde{x}\ll 1.
\label{eq:skellam}
\end{equation}
It is natural to remove linear source term $\lambda\tilde{u}$ in (\ref{eq:skellam}) by defining a new independent variable
\begin{equation}
\tilde{v}(\tilde{x},\tilde{t})=\mathrm{e}^{\lambda(\tilde{t}_{\mathrm{c}}-\tilde{t})}\,\tilde{u},
\label{eq:newindep}
\end{equation}
so that $\tilde{v}$ satisfies the heat equation subject to
\begin{equation}
\tilde{v}=0, \quad \frac{\mathrm{d}\tilde{s}}{\mathrm{d}\tilde{t}}=-\kappa\,
\mathrm{e}^{-\lambda(\tilde{t}_{\mathrm{c}}-\tilde{t})}
\frac{\partial \tilde{v}}{\partial \tilde{x}}
\quad\mbox{on}\quad \tilde{x}=\tilde{s}(\tilde{t}).
\label{eq:newBCs}
\end{equation}
In this way, we can see the close connection between our moving problem with the Fisher-KPP equation and the Stefan problem for the heat equation (see section~\ref{sec:superheated}).

Indeed, for $\lambda$ sufficiently small, we follow the ideas by King \& Evans~\cite{King2005} for the Stefan problem although, as mentioned, we shall not use the Baiocchi transform.  First, we look for a solution of the form
$$
\tilde{v}\sim \frac{1}{(\tilde{t}_{\mathrm{c}}-\tilde{t})^a}f(\xi,\tau)
\quad\mbox{as}\quad \tilde{t}\rightarrow \tilde{t}_{\mathrm{c}}^-,
$$
where the similarity variables are
$$
\xi=\frac{\tilde{x}-\tilde{s}_{\mathrm{c}}}{(\tilde{t}_{\mathrm{c}}-\tilde{t})^b},\quad \tau=-\ln(\tilde{t}_{\mathrm{c}}-\tilde{t}).
$$
Immediately, satisfying the heat equation for $\tilde{v}$ means that $b=1/2$.  The choice $a=0$ allows a match in the boundary conditions (\ref{eq:StefanFisherscaled3}).  The result is that
\begin{equation}
\tilde{v}\sim f(\xi,\tau)
\quad\mbox{as}\quad \tau\rightarrow\infty,
\label{eq:vf}
\end{equation}
where
\begin{equation}
\xi=\frac{\tilde{x}-\tilde{s}_{\mathrm{c}}}{(\tilde{t}_{\mathrm{c}}-\tilde{t})^{1/2}},\quad \tau=-\ln(\tilde{t}_{\mathrm{c}}-\tilde{t}),
\quad \sigma(\tau)=\frac{\tilde{s}(\tilde{t})-\tilde{s}_{\mathrm{c}}}{(\tilde{t}_{\mathrm{c}}-\tilde{t})^{1/2}},
\label{eq:similarity}
\end{equation}
and $f$ satisfies
\begin{equation}
\frac{\partial f}{\partial\tau}=\frac{\partial^2 f}{\partial \xi^2}-\frac{1}{2}\xi\frac{\partial f}{\partial \xi}, \quad \xi<\sigma(\tau),
\label{eq:fpde}
\end{equation}
\begin{equation}
f=0,\quad \dot{\sigma}-\frac{1}{2}\sigma = -\kappa \,\mathrm{e}^{-\lambda\mathrm{exp}(-\tau)}f_\xi
\quad\mbox{on}\quad \xi=\sigma(\tau),
\label{eq:fBCs}
\end{equation}
where here the dot means a derivative with respect to $\tau$.  The use of self-similar variables $\xi$ and $\tau$ in (\ref{eq:vf})-(\ref{eq:similarity}) (instead of just $\xi$) is discussed at length in Ref.~\cite{Eggers2009}, for example.

\begin{figure}
	\centering
\begin{tikzpicture}[xscale=1 ,yscale=1,cap=round]
font=\footnotesize,
thin,
\draw (0,0) -- (2,0);

\draw (1.75,-0.25) -- (2.25,0.25);
\draw (2.25,-0.25) -- (2.75,0.25);

\draw[->,-latex] (2.5,0) -- (11.5,0) coordinate(x axis) node[right] {$\xi$}; 

\draw (0,0) -- (0,-0.1) node[below] {$ \displaystyle \lim_{\tilde{t}\to \tilde{t}_\mathrm{c}^-}\frac{-\tilde{s}_\mathrm{c}}
{(\tilde{t}_\mathrm{c}-\tilde{t})^{1/2}} = -\infty$};

\draw (11,0) -- (11,-0.1);
\draw (10.5,-0.1) node[below] {$ \displaystyle \lim_{\tilde{t}\to \tilde{t}_\mathrm{c}^-}\sigma(\tau) = \infty$};

\draw (4.75,0) -- (4.75,-0.1) node[below] {$0$};
\draw[dashed] (3.75,0.1) -- (3.75,3);
\draw[dashed] (5.75,0.1) -- (5.75,3);
\draw[dashed] (7.5,0.1) -- (7.5,3);
\draw[dashed] (11,0.1) -- (11,3);

\draw (1.5,1.5) node {$\xi < 0$};
\draw (4.75,2) node {outer};
\draw (4.75,1.5) node {$\xi = \mathcal{O}(1)$};
\draw (6.6,2) node {matching};
\draw (6.75,1.5) node {$\xi \gg 1$};
\draw (9.25,2.5) node {interior layer};
\draw (9.25,1.75) node {$\displaystyle \xi = \sigma(\tau) + \frac{\eta}{\sigma(\tau)}$,};
\draw (9.25,1) node {$\eta= \mathcal{O}(1)$};
\end{tikzpicture}
	\caption{A schematic of the structure of the near-blow-up analysis in terms of the variables $\xi=(\tilde{x}-\tilde{s}_{\mathrm{c}})/(\tilde{t}_{\mathrm{c}}-\tilde{t})^{1/2}$ and $\tau=-\ln(\tilde{t}_{\mathrm{c}}-\tilde{t})$.  The diminishingly small region $\tilde{s}_{\mathrm{c}}<\tilde{x}<\tilde{s}(\tilde{t})$ is stretched out to $0<\xi<\sigma(\tau)$, where $\sigma\rightarrow\infty$ as $\tilde{t}\rightarrow\tilde{t}_{\mathrm{c}}^-$ ($\tau\rightarrow\infty$).}
	\label{fig:schematic}
\end{figure}
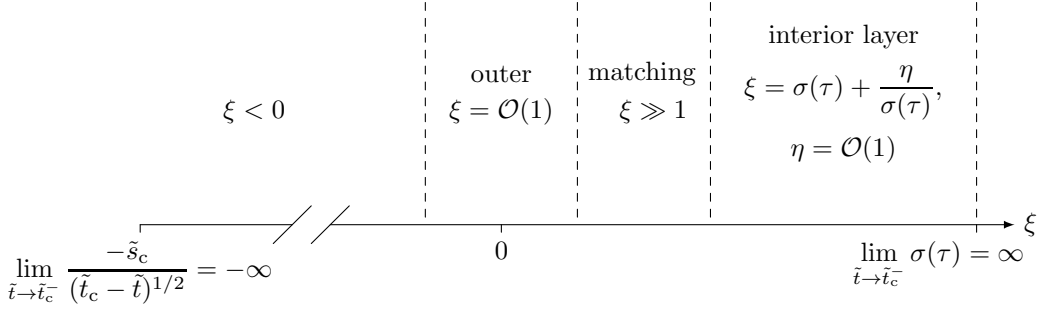

Before continuing, it is worth explaining that a key reason the analysis for this problem is much more complicated than a traditional similarity solution is that the straightforward ansatz $\tilde{v}\sim f(\xi)$ as $\tau\rightarrow\infty$ does not work, as it leads to $f''-\xi f'/2=0$, whose solutions blow up exponentially fast as $\xi\rightarrow -\infty$.  The culprit is the sign of $\tilde{t}$ in $(\tilde{t}_{\mathrm{c}}-\tilde{t})^{1/2}$, which leads to an unwanted negative sign when taking the derivative with respect to time.  Further, for such a straightforward similarity solution, the moving boundary would be characterised by $\sigma=\,\mbox{constant}$, giving $\tilde{s}-\tilde{s}_{\mathrm{c}}\sim \mbox{constant}\, (\tilde{t}_{\mathrm{c}}-\tilde{t})^{1/2}$; however, as stated above in (\ref{eq:blowup}), it turns out that $\sigma\rightarrow\infty$ (very slowly) in the limit, leading to a complicated spatial structure for $\xi$ in the limit $\tau\rightarrow\infty$, as we now summarise.  A schematic to help with this description is provided in figure~\ref{fig:schematic}.

Following the key ideas in Ref.~\cite{King2005}, for $\xi=\mathcal{O}(1)$ we write
\begin{equation}
f\sim C + A(\tau)+\dot{A}(\tau)\Phi_0(\xi)+\ddot{A}(\tau)\Phi_1(\xi)
\quad\mbox{as}\quad \tau\rightarrow\infty,
\label{eq:outer}
\end{equation}
where $C$ is a constant and $A(\tau)$ is an unknown function with the property $|\ddot{A}|\ll |\dot{A}|\ll |A| \ll 1$ as $\tau\rightarrow\infty$.  By substituting (\ref{eq:outer}) into (\ref{eq:fpde}), to leading order we find
$$
\Phi_0''-\frac{1}{2}\xi\Phi_0'=1.
$$
Excluding exponential growth as $\xi\rightarrow -\infty$, we find $\Phi_0'=\sqrt{\pi}\mathrm{e}^{\xi^2/4}\left(1+\mathrm{erf}(\xi/2)\right)$, with the asymptotic behaviour $\Phi_0'\sim -2/\xi$ as $\xi\rightarrow -\infty$ and $\Phi_0'\sim 2\sqrt{\pi}\mathrm{e}^{\xi^2/4}-2/\xi+4/\xi^3$ as $\xi\rightarrow\infty$. Therefore we have, as a matching condition,
\begin{equation}
\Phi_0\sim \frac{4\sqrt{\pi}}{\xi}\,\mathrm{e}^{\xi^2/4}\quad\mbox{as}\quad \xi\rightarrow\infty.
\label{eq:largexi}
\end{equation}
The solution for $\Phi_1$ is not required in what follows.  Note that expansions like (\ref{eq:outer}), together with an ansatz like (\ref{eq:vf})-(\ref{eq:similarity}), are used in studies of singularities in other Stefan problems, for example for inward solidification \cite{Herrero1997,McCue2005,Soward1980}.

Near the moving boundary $\xi=\sigma(\tau)$ there is an interior layer with
$$
\xi=\sigma+\frac{\eta}{\sigma},
$$
where $\eta=\mathcal{O}(1)$.  In this boundary layer, let $f=\Psi(\eta,\tau)$, so that (\ref{eq:fpde}) becomes
\begin{equation}
\frac{1}{\sigma^2}\frac{\partial\Psi}{\partial\tau}+\left(-\sigma+\frac{\eta}{\sigma}\right) \frac{\dot{\sigma}}{\sigma^2}\frac{\partial\Psi}{\partial\eta}
=\frac{\partial^2\Psi}{\partial\eta^2}-\frac{1}{2}\left(1+\frac{\eta}{\sigma^2}\right)
\frac{\partial\Psi}{\partial\eta},
\quad \eta<0,
\label{eq:interiorpde}
\end{equation}
which is dominated by the right-hand side in the limit $\tau\rightarrow\infty$.  We write
\begin{equation}
f\sim \Psi_0(\eta)+\frac{1}{\sigma^2}\Psi_1(\eta)
\quad\mbox{as}\quad \tau\rightarrow\infty,
\label{eq:interior}
\end{equation}
so that, after substituting into (\ref{eq:interiorpde}), we find
$$
\Psi_0''-\frac{1}{2}\Psi_0'=0,\quad
\Psi_1''-\frac{1}{2}\Psi_1'=\frac{1}{2}\eta\Psi_0'.
$$
After satisfying the first condition in (\ref{eq:fBCs}), we find $\Psi_0=\mathrm{constant}(1-\mathrm{e}^{\eta/2})$.  To determine the constant, we assume that
\begin{equation}
\lambda\,\mathrm{exp}(-\tau)\ll 1,
\quad\mbox{or}\quad
\tilde{t}_{\mathrm{c}}-\tilde{t}\ll \lambda^{-1},
\label{eq:lambdatau}
\end{equation}
so that the second condition in (\ref{eq:fBCs}) gives $-1/2=-\kappa\Psi_0'(0)$, which means
\begin{equation}
\Psi_0=-\frac{1}{\kappa}(1-\mathrm{e}^{\eta/2}).
\label{eq:Psi0}
\end{equation}
The solution for $\Psi_1$ subject to $\Psi_1=0$ and $\Psi_1'=0$ at $\eta=0$ (from the first and second conditions in (\ref{eq:fBCs}), respectively) is
\begin{equation}
\Psi_1=-\frac{2}{\kappa}+\frac{\eta^2-4\eta+8}{4\kappa}\,\mathrm{e}^{\eta/2}.
\label{eq:Psi1}
\end{equation}

Matching between the interior layer (\ref{eq:interior}), (\ref{eq:Psi0})-(\ref{eq:Psi1}) and the region for $\xi=\mathcal{O}(1)$, (\ref{eq:outer}), in the limit $\eta\rightarrow -\infty$, gives
$$
C=-\frac{1}{\kappa},\quad A\sim -\frac{2}{\kappa \sigma^2},
$$
while matching between the exponential terms in (\ref{eq:largexi}) and (\ref{eq:Psi0}) gives
$$
\frac{4\sqrt{\pi}}{\sigma}\,\mathrm{e}^{\sigma^2/4}\dot{A}\sim \frac{1}{\kappa}.
$$
Together, these results give
$$
\frac{16\sqrt{\pi}}{\kappa\sigma^4}\,\mathrm{e}^{\sigma^2/4}\dot{\sigma}\sim \frac{1}{\kappa},
\quad\mbox{or}\quad
\frac{32\sqrt{\pi}}{\sigma^5}\,\mathrm{e}^{\sigma^2/4}\sim \tau,
\quad\mbox{as}\quad \tau\rightarrow\infty.
$$
Finally, we find that, to leading order,
$$
A\sim -\frac{1}{2\kappa\ln\tau},
\quad
\sigma\sim 2\ln^{1/2}\tau,
\quad\mbox{as}\quad \tau\rightarrow\infty.
$$
The latter implies that $\tilde{s}(\tilde{t})-\tilde{s}_{\mathrm{c}}\sim 2(\tilde{t}_{\mathrm{c}}-\tilde{t})^{1/2}\ln^{1/2}\tau$, which is (\ref{eq:blowup}).

In the original variables, the blow-up profile has the limiting behaviour
\begin{equation}
\tilde{u}(\tilde{x},\tilde{t}_{\mathrm{c}})\sim -\frac{1}{\kappa}\left(1+
\frac{1}{2\ln\left(-\ln(\tilde{s}_{\mathrm{c}}-\tilde{x})
\right)}\right)
\quad\mbox{as}\quad \tilde{x}\rightarrow \tilde{s}_{\mathrm{c}}^-.
\label{eq:blowupprofile}
\end{equation}
This is the same limiting profile as for the Stefan problem described in subsection~\ref{sec:superheated}.  As mentioned above, the key reason that (\ref{eq:blowup}) and (\ref{eq:blowupprofile}) are unchanged for $\lambda=\mathcal{O}(1)$ when compared to $\lambda=0$ is that the reaction terms $\lambda\tilde{u}(1-\tilde{u})$ do not contribute to leading order behaviour near the moving boundary $\tilde{x}=\tilde{s}(\tilde{t})$.

\begin{figure}
	\centering
	\includegraphics[width=1\linewidth]{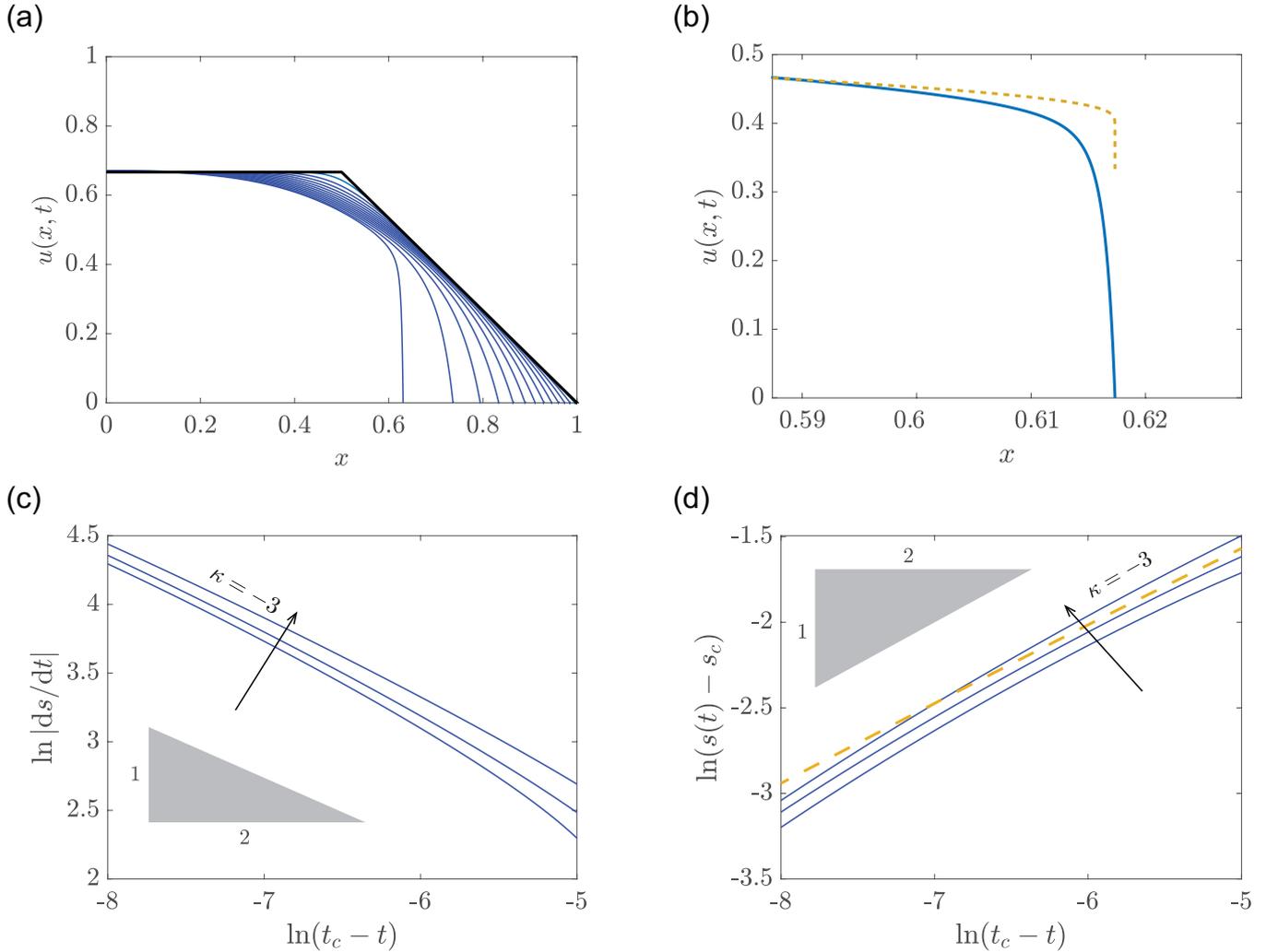}
	\caption{ (a) Time-dependent solutions of (\ref{eq:StefanFisherscaled1})-(\ref{eq:StefanFisherscaled4}) (in blue) for $\kappa = -3$, $\lambda=1$ and the ramp style initial condition with $M(0)=0.5$ (in black).  The solutions are drawn at time intervals of $2\times10^{-3}$. (b) Solution at $\tilde{t}_{\mathrm{c}}\approx 0.024$ (in blue) is compared to (\ref{eq:blowupprofile}) (orange dashed) where $\tilde{s}_c=0.614$.  Log-log plots for (c) $\mathrm{d}\tilde{s}/\mathrm{d}\tilde{t}$ versus $\tilde{t}_{\mathrm{c}}-\tilde{t}$ and (d) $\tilde{s}(\tilde{t})-\tilde{s}_{\mathrm{c}}$ versus $\tilde{t}_{\mathrm{c}}-\tilde{t}$, drawn for $\kappa=-5$, $-4$, $-3$. Arrows indicate increasing values of $\kappa$. Also included is a line segment with slope $-1/2$ in (c) and $1/2$ in (d).  The near blow-up result (\ref{eq:blowup}) is shown in (d) as an orange dashed curve.}
	\label{figure8}
\end{figure}

Limits like (\ref{eq:blowup}) and (\ref{eq:blowupprofile}) are incredibly difficult to test numerically, as are all blow-up scenarios where the speed of the moving boundary becomes unbounded in a finite time.  In this case, it is the double logarithm that grows so slowly that its effect kicks in on time and length scales that are pushing the limits of the numerical scheme.  In figure~\ref{figure8} we show simulations of (\ref{eq:StefanFisherscaled1})-(\ref{eq:StefanFisherscaled4}) for $\lambda=1$ and values of $\kappa$ that lead to blow-up.  In figure~\ref{figure8}(a), solutions are shown for various times for $\tilde{t}>0$ (in blue), including the initial condition (in black), which is a ramp style function.  In figure~\ref{figure8}(b), part of a solution profile is drawn close to the estimated blow-up time, again as a solid blue curve, together with (\ref{eq:blowupprofile}), which is a thick orange dashed curve.  Given the challenges involved in making this comparison, the agreement is quite good.

In terms of the time-dependence of solutions shown in figure~\ref{figure8}(a)-(b), we have included log-log plots of both $\mathrm{d}\tilde{s}/\mathrm{d}\tilde{t}$ and $\tilde{s}(\tilde{t})-\tilde{s}_{\mathrm{c}}$ versus $\tilde{t}_{\mathrm{c}}-\tilde{t}$ in figure~\ref{figure8}(c) and (d), respectively, computed for three different values of $\kappa$.  Included in these images are line segments with slopes $-1/2$ and $1/2$, to indicate the dominant part of the power-law dependence in (\ref{eq:blowup}), together with (\ref{eq:blowup}) itself.  Of course it is difficult to capture this behaviour on very small time-scales, especially given the logarithm of the logarithm in (\ref{eq:blowup}), which only begins to grow significantly for extremely small values of $\tilde{t}_{\mathrm{c}}-\tilde{t}$.  With this in mind, the agreement is acceptable.

\section{Discussion}\label{sec:conclusion}

In this paper we have studied a Stefan-type moving boundary problem where the usual heat equation is replaced by the Fisher-KPP equation.  We refer to this system as the Fisher-Stefan model.  From a biological perspective, modifying the Fisher-KPP equation to include a moving boundary has a number of advantages, for example it allows for solutions with compact support which are more realistic than those that decay exponentially in the far field.  This particular advantage is also shared with models that use nonlinear degenerate diffusion for cell motility~\cite{Simpson2011}, for example resulting in the Porous-Fisher equation \cite{Fadai2020,McCue2019,Sanchez1995,Witelski1995}, without needing to discard linear diffusion.  Various other moving boundary problems that are adapted from the Fisher-KPP model have been studied in recent times~\cite{Berestycki2018,Berestycki2019,Berestycki2021,Tisbury2021a}.

The Fisher-Stefan model (\ref{eq:StefanFisher1})-(\ref{eq:StefanFisher4}) involves a new parameter $\kappa$, which is a measure of leakage at the moving boundary.  This extra parameter allows for a range of qualitative behaviours that are not available for the usual Fisher-KPP model.  In particular, for $\kappa>0$, the moving boundary invades an empty space and, depending on the initial conditions, the solution may evolve to travelling wave or go to extinction.  Du \& Lin~\cite{Du2010} refer to this situation as the spreading-vanishing dichotomy.  Aspects of these phenomena have been studied rigorously and numerically by Du and coworkers \cite{Du2010,Du2020,Du2014a,Du2014b,Du2015,Liu2020}, us \cite{Elhachem2019,Elhachem2020,Elhachem2021,McCue2020,Simpson2020} and others.

Our study in this paper is for $\kappa<0$, for which the moving boundary retreats from the empty space (see also \cite{Elhachem2021,McCue2020}).  By neglecting the effect of imposing a boundary condition at the left-hand boundary (section~\ref{sec:blowup}), we have used numerical simulation to demonstrate how solutions for $-1<\kappa<0$ evolve to retreating travelling waves while those for $\kappa<-1$ blow up in finite time; for the borderline case $\kappa=-1$, it appears that solutions blow up in infinite time.  By including the effects of a no-flux condition at the left boundary (section~\ref{sec:stefan}), we eliminate the possibility of travelling wave solutions.  Instead, the solution may also exhibit finite-time blow-up for $\kappa$ less than some critical value, but this critical value depends on the initial condition and the proliferation rate $\lambda$.  Further, the other option apart from blow-up is extinction, which occurs when the initial population size (or ``mass'') is sufficiently small.  There is a direct analogy here with the problem of melting a superheated solid \cite{Howison1985,Back2010,Fasano1977,Fasano1981,Herrero1996,King2005,Lacey1985,Sherman1970}, which we describe in some detail.

The finite-time blow-up we report on here involves the speed of the moving boundary and the slope of the solution at the moving boundary both becoming infinite at some time $t_\mathrm{c}<\infty$.  The population density itself remains bounded in the interval $0\leq u\leq 1$.  One of the reasons why this study is interesting is that blow-up in reaction diffusion equations is more commonly accompanied by the population itself increasing with outbound, for example where the reaction term is of the form $u^p$ where $p>1$~\cite{Deng2000,Levine1990}.  Blow-up with a bounded dependent variable, as we observe here with a model based on a second-order parabolic partial differential equation, is more likely to be observed in first-order wave equations where shocks are well-known to form~\cite{Whitham2011}.

We have presented a summary of the blow-up asymptotics for the rescaled problem (\ref{eq:StefanFisherscaled1})-(\ref{eq:StefanFisherscaled4}) in section~\ref{sec:blowuprevisited}.  After ignoring the nonlinear term $-\tilde{u}^2$, the details follow that for the superheated Stefan problem with $\lambda=0$ closely \cite{Herrero1996,King2005}, although we have presented the analysis without resorting to the Baiocchi transform.  A key observation is that the time-dependence (\ref{eq:blowup}) relies on (\ref{eq:lambdatau}), which means that if $\lambda=\mathcal{O}(1)$, then (\ref{eq:blowup}) provides a good leading order estimate.  However, if $\lambda\gg 1$, then (\ref{eq:blowup}) holds for such an incredibly small time interval before blow-up that it is too difficult to test numerically.  This is an important observation because, due to (\ref{eq:lambda}), our numerical simulations in figures~\ref{figure1}, \ref{figure3}-\ref{figure4} are for $\lambda=10^6$, which is very large.

In terms of further work, it would be interesting to revisit our numerical results in a rigorous setting.  For example, by posing (\ref{eq:StefanFisher1})-(\ref{eq:StefanFisher4}) on $-\infty<x\leq s(t)$,
a proof of global existence for $-1\leq\kappa<0$ and local existence for $\kappa<-1$ as worth exploring,
as is a rigorous exploration of the near blow-up conjecture (\ref{eq:blowup}).  Further, a proof of our conjecture that the borderline case $\kappa=-1$ corresponds to infinite-time blow-up (on the semi-infinite domain) would be welcome.  In terms of extinction, the rough estimate (\ref{eq:prection2}) could be refined using careful bounds on the solution with $\lambda>0$.

Another issue worth further formal investigation is stability.  In the Stefan problem (\ref{eq:StefanFisherscaled1})-(\ref{eq:StefanFisherscaled4}) with $\lambda=0$, when $\kappa>0$ (normal melting/freezing) the solutions are stable to perturbations in the transverse direction (in the $y$-direction, say), while for $\kappa<0$ (melting a superheated solid or freezing a supercooled liquid) the solutions are highly unstable, even with a surface-tension-type regularisation \cite{Mullins1963}.  It would be interesting from a mathematical perspective to extend known stability results to the Fisher-Stefan model with a moving boundary.  Numerical simulations and linear stability analysis would provide insight whether the front is unstable and which modes grow the fastest, for example.  Such an approach would require an analogous restoring force like surface tension, as included in other moving boundary problems in mathematical biology~\cite{Giverso2015}.

\paragraph{Acknowledgements:} This work is supported by the Australian Research Council (DP200100177).  SWM is grateful to John King and Julian Back for discussions on Stefan problems.  The authors thank the anonymous referees for their constructive feedback.

\appendix

\section{Numerical methods} \label{sec:Numericalmethods}

To obtain numerical solutions of the Fisher--Stefan equation
\begin{equation}\label{eq:FisherKPPmovbound}
\frac{\partial u}{\partial t} =\frac{\partial^2 u}{\partial x^2} + \lambda u(1-u),
\end{equation}
for $0 < x < s(t)$ and $t > 0$, we use a boundary fixing transformation $y = x / s(t)$ so that we have
\begin{align}\label{eq:FisherKPPmovboundxi}
\frac{\partial u}{\partial t} = \frac{1}{s^2(t)} \frac{\partial^2 u}{\partial y^2}+\frac{y}{s(t)} \frac{\text{d}s(t)}{\text{d} t} \frac{\partial u}{\partial y} + \lambda u(1-u),
\end{align}
on the fixed domain, $0 < y < 1$.  Here $s(t)$ is the time--dependent length of the domain, and we will explain how we solve for this quantity later. To close the problem we also transform the boundary conditions giving
\begin{align}
&\dfrac{\partial u}{\partial y} = 0 \quad \textrm{at} \quad y=0,  \label{eq:FS_BC1} \\
&u = 0 \quad \textrm{at} \quad y=1. \label{eq:FS_BC2}
\end{align}
The key to obtaining accurate numerical solutions of equation (\ref{eq:FisherKPPmovbound}) is to take advantage of the fact that for many problems we consider $u(x,t)$ varies rapidly near $x=s(t)$, whereas $u(x,t)$ is approximately constant near $x=0$.  Motivated by this spatial structure, we discretise equation (\ref{eq:FisherKPPmovboundxi}) using a variable mesh where the mesh spacing varies geometrically from $\delta y_{\textrm{min}} = y_{N} - y_{N-1} = 1 - y_{N-1}$ at $y = 1$, to $\delta y_{\textrm{max}} = y_2 - y_1 = y_2 - 0$ at $y = 0$.  Results in this work are computed with $N=1001$ mesh points with $\delta y_{\textrm{min}} = 1 \times 10^{-6}$.  With these constraints we solve for the geometric expansion factor using MATLABs \textit{fsolve} function which gives $\delta y_{\textrm{max}} = 9.083 \times 10^{-2}$ and the geometric expansion factor is 1.009165 to six decimal places.

We spatially discretise equation (\ref{eq:FisherKPPmovboundxi}) on the non-uniform mesh.  At the $i$th internal mesh point we define $h_i^+ = y_{i+1} - y_{i}$ and $h_i^- = y_{i} - y_{i-1}$.   For convenience we define $\alpha_i = 1/(h^-[h^+ + h^-])$, $\gamma_i = -1/(h^- h^+)$ and $\delta_i = 1/(h^+[h^+ + h^-])$, which gives
\begin{align}\label{eq:FDDinternalxi}
\dfrac{u_{i}^{j+1} - u_{i}^{j} }{\Delta t} &=  \dfrac{2}{(s^{j})^2} \left[ \alpha_i u_{i-1}^{j+1} + \gamma_i u_{i}^{j+1} + \delta_i u_{i+1}^{j+1}  \right]\notag \\
&+\dfrac{y_i}{s^{j}} \left(\dfrac{s^{j+1} - s^{j} }{\Delta t}\right) \left[ \delta_i h^N u_{i-1}^{j+1} + \gamma_i (h^N-h^+) u_{i-1}^{j+1} - \alpha_i h^+ u_{i+1}^{j+1} \right] \notag \\
& +  \lambda u_{i}^{j+1}(1 - u_{i}^{j+1}) ,
\end{align}
for $i = 2, \ldots, N-1$, where $N$ is the total number of spatial nodes on the finite difference mesh, and the index $j$ represents the time index so that $u_{i}^{j} \approx u(y_i, j\Delta t)$.
Finally, discretising equations (\ref{eq:FS_BC1})--(\ref{eq:FS_BC2}) leads to
\begin{align}
&u_{2}^{j+1}-u_{1}^{j+1} = 0,  \label{eq:FS_BC1a} \\
&u_{N}^{j+1} = 0. \label{eq:FS_BC2a}
\end{align}

To advance the discrete system from time $t$ to $t + \Delta t$ we solve the system of nonlinear algebraic equations, equations (\ref{eq:FDDinternalxi})-(\ref{eq:FS_BC2a}), using Newton-Raphson iteration. During each iteration of the Newton--Raphson algorithm we estimate the position of the moving boundary using the discretised Stefan condition.  Here we define $h_N^+ = y_{N} - y_{N-1}$, $h_N^- = y_{N-1} - y_{N-2}$, $\alpha_N = 1/(h^-[h^+ + h^-])$ and $\gamma_N = -1/(h^- h^+)$, which gives
\begin{equation}
s^{j+1} = s^{j} - \dfrac{\Delta t  \kappa}{s^j} \left[\alpha_N h^+ u_{N-1}^{j+1} + \gamma_N (h^+ + h^-) u_{N}^{j+1} \right].
\label{eq:Lupdatediscretise}
\nonumber
\end{equation}
Within each time step Newton--Raphson iterations continue until the maximum change in the dependent variables is less than the tolerance  $\epsilon$.  All results in this work are obtained by setting $\epsilon = 1 \times 10^{-10}$, and $\Delta t = 1 \times 10^{-4}$, and we find that these values are sufficient to produce grid--independent results.  Software written in MATLAB is available on \href{https://github.com/maudelhachem/El-Hachem2020b}{GitHub} for the reader to experiment with different choices of discretization and tolerance values.

\end{document}